\documentstyle[amstex]{article}
\begin{document}

\title{{\bf   Realcompactness and spaces of vector-valued functions}}

\author{\normalsize {\sc Jes\'us Araujo}\thanks{Research partially 
supported by the
Spanish Direcci\'on General de Investigaci\'on Cient\'{\i}fica y
T\'ecnica (DGICYT,  PB98-1102).}}
% {\sc \hspace{0.04in} and \hspace{0.04in} 
%Juan J. Font} \thanks{Research partially supported by Spanish DGES grant 
%(PB96-1075).}}

\date{}                                         

\maketitle 

%\newtheorem{lemma}{Lemma}
%\newtheorem{theorem}{Theorem}
%\newtheorem{corollary}{Corollary}
%\newtheorem{proposition}{Proposition}
%\newtheorem{example}{Example}
%\newtheorem{remark}{Remark}

%\documentstyle[12pt]{article}

%\newcommand{\ntt}{\series m\shape n\tt}

%\newcommand{\cs}[1]{{\protect\ntt\bslash#1}}

%\newcommand{\opt}[1]{{\protect\ntt#1}}
%\newcommand{\env}[1]{{\protect\ntt#1}}

%\makeatletter
%\def\verbatim{\interlinepenalty\@M \@verbatim
%  \leftskip\@totalleftmargin\advance\leftskip2pc
%  \frenchspacing\@vobeyspaces \@xverbatim}
%\makeatother
%\hfuzz1pc 

%       Theorem environments

\newtheorem{theorem}{Theorem}[section]
\newtheorem{corollary}[theorem]{Corollary}
\newtheorem{lemma}[theorem]{Lemma}
\newtheorem{prop}[theorem]{Proposition}
\newtheorem{ax}{Axiom}

\newtheorem{definition}{Definition}[section]

\numberwithin{equation}{section}

\newtheorem{examples}{Examples}

\newtheorem{example}{Example}
%\renewcommand{\theexample}{}

%       Math definitions

\newcommand{\cly}{{\rm cl}_{\beta Y} \hspace{.02in}}
\newcommand{\cgy}{{\rm cl}_{\gamma Y} \hspace{.02in}}
\newcommand{\chy}{{\rm cl}_{Y} \hspace{.02in}}
\newcommand{\clx}{{\rm cl}_{\beta X} \hspace{.02in}}
\newcommand{\cgx}{{\rm cl}_{\gamma X} \hspace{.02in}}
\newcommand{\clp}{{\rm cl}_{{\Bbb R}^p} \hspace{.02in}}
\newcommand{\chx}{{\rm cl}_{X} \hspace{.02in}}
\newcommand{\cux}{{\rm cl}_{\upsilon X} \hspace{.02in}}
\newcommand{\bx}{A(X,E)}
\newcommand{\xe}{C^*(X,E)}
\newcommand{\yf}{C^*(Y,F)}
\newcommand{\ra}{\rightarrow}
\newcommand{\smn}{\sum_{n=1}}
\newcommand{\ay}{A(Y,F)}
\newcommand{\cxo}{C_0(X,E)}
\newcommand{\cyo}{C_0(Y,F)}
\newcommand{\cx}{C_0(X)}
\newcommand{\sxy}{T:A(X) \ra A(Y)}
\newcommand{\txy}{T:A(X,E) \ra A(Y,F)}
\newcommand{\bxy}{T:C^* (X,E) \ra C^* (Y,F)}
\newcommand{\uxy}{T:C(X,E) \ra C(Y,F)}
\newcommand{\vxy}{T:C_0 (X,E) \ra C_0 (Y,F)}
\newcommand{\pa}{\left(}
\newcommand{\rb}{\right)}
\newcommand{\va}{\left|}
\newcommand{\vb}{\right|}
\newcommand{\vc}{\left\|}
\newcommand{\vd}{\right\|}
\newcommand{\cl}{{\rm cl} \hspace{.02in}}
\newcommand{\itr}{{\rm int} \hspace{.02in}}
\newcommand{\itrx}{{\rm int}_{\beta X} \hspace{.02in}}
\newcommand{\bax}{^{\beta X}}
\newcommand{\bay}{^{\beta Y}}
\newcommand{\ing}{{\rm ifn} \hspace{.02in}}
\newcommand{\fin}{{\rm fin} \hspace{.02in}}
\newcommand{\bdry}{{\rm bdry} \hspace{.02in}}
\newcommand{\hs}{\hspace{.02in}}

\newcommand{\noe}{A^n (\Omega, E)}
\newcommand{\coe}{C^n (\Omega, E)}
\newcommand{\no}{C^n (\Omega, {\Bbb  R})}
\newcommand{\nor}{C^n ({\Bbb  R}, {\Bbb  R})}
\newcommand{\noor}{C^n_c (\Omega, {\Bbb  R})}
\newcommand{\mor}{C^m ({\Bbb  R}, {\Bbb  R})}
\newcommand{\moor}{C^m_c (\Omega', {\Bbb  R})}
\newcommand{\mo}{C^m (\Omega', {\Bbb  R})}
\newcommand{\moe}{A^m (\Omega', F)}
\newcommand{\bnoe}{C^n_* (\Omega, E)}
\newcommand{\bno}{A^n (\Omega, {\Bbb  R})}
\newcommand{\bmo}{C^m_* (\Omega', {\Bbb  R})}
\newcommand{\bmoe}{C^m_* (\Omega', F)}
\newcommand{\ixy}{T: \bnoe \ra \bmoe}
\newcommand{\partl}{\frac{\partial^{\va \lambda \vb} l}{\partial x^{\lambda}}}
\newcommand{\parto}{\frac{\partial^{\va \lambda \vb} l'}{\partial x^{\lambda}}}
\newcommand{\partm}{\frac{\partial^{\va \lambda \vb} l}{\partial x^{\lambda}}}
\newcommand{\partg}{\frac{\partial^{\va \lambda \vb} Tg}{\partial x^{\lambda}}}
\newcommand{\partf}{\frac{\partial^{\va \lambda \vb} Tf}{\partial x^{\lambda}}}

%\title[Banach-Stone theorems]{Realcompactness and Banach-Stone theorems}

%\thispagestyle{empty}

%% First author
%\author{Jes\'us Araujo}
%\address{Departamento de Matem\'aticas,
%Estad\'{\i}stica y Computaci\'on\\ Universidad de Cantabria\\
%Facultad de Ciencias\\ Avda.
%de los Castros, s. n.\\ E-39071 Santander, Spain}
%% Note the doubled @@:
%\email{araujo@@matesco.unican.es}
%\thanks{Research partially supported by
%the Spanish Direcci\'on General de Investigaci\'on Cient\'{\i}fica
%y T\'ecnica (DGICYT, PB98-1102).}
%\thanks{Research of the first author was %supported in part by
%the Spanish Direcci\'on General de Investigaci\'on Cient\'{\i}fica
%y T\'ecnica (DGICYT, PS90-100).}

%\subjclass{Primary 46E40; Secondary 47B33, 47B38, 54D60}
%\keywords{Biseparating map, Stone-\v{C}ech compactification, realcompactification, %linear isometry, Banach-Stone theorem}

\date{}                                          
\thispagestyle{empty}
%\begin{document}

\maketitle

\begin{abstract}
It is shown that the existence of a biseparating map between
a large class of spaces of vector-valued continuous functions $\bx$ and
$\ay$ implies that  some compactifications of $X$ and $Y$ are
homeomorphic. In some cases, conditions are given to warrant the existence of a homeomorphism between the  realcompactifications of $X$ and $Y$; in particular we find remarkable differences with respect to the scalar context: namely, if $E$ and $F$ are infinite-dimensional and $T: C^{*} (X,E) \ra C^{*} (Y, F)$ is a biseparating map, then the realcompactifications of $X$ and $Y$ are homeomorphic.
\end{abstract}

\section{Introduction}
Let ${\Bbb K} = {\Bbb R}$ or ${\Bbb C}$.
Given a completely regular space $X$, and a ${\Bbb K}$-normed space $E$,
$C(X,E)$ and $\xe$ denote the spaces of continuous
functions and {\em
bounded} continuous functions on $X$ taking 
values on $E$, respectively. $C(X)$ and $C^* (X) $ will be the
spaces $C(X, {\Bbb K})$ and $C^* (X, {\Bbb K})$, respectively.

        Sometimes an algebraic relation between spaces of continuous
functions $C(X)$ and $C(Y)$ may determine some kind of
topological link between the spaces $X$ and $Y$. For
instance, it is well known that the existence of a ring
isomorphism between the spaces $C(X)$ and $C(Y)$ produces a
homeomorphism between the realcompactifications of $X$ and $Y$
(see \cite[pp. 115-118]{GJ} and \cite{H}). Some kind of weakening on the
conditions do not alter the result: if we replace "ring
isomorphism" by "biseparating map" (see definition below), we
keep the conclusion on the existence of a homeomorphism between the realcompactifications of
$X$ and $Y$ (\cite{ABN1}). On the other hand, it is also true that if the realcompactifications $\upsilon X$ of $X$ and $\upsilon Y$ of $Y$ are homeomorphic through a homeomorphism $h:\upsilon Y \ra \upsilon X$, then there exists a biseparating map (in fact, a ring isomorphism) between $C(X)$ and $C(Y)$, as the map sending each $f \in C(\upsilon X)$ into $f \circ h \in C(\upsilon Y)$ can be identified with an isomorphism between $C(X)$ and $C(Y)$.
%Also, in
%some cases, if we assume linearity on $T: C(X) \ra C(Y)$, it is possible to
%wither even more the conditions on $T$ without getting any
%change on the conclusion (\cite{AK, ABN2, J}). In general, linearity
%of $T$ gives even a more interesting property, as it is the
%possibility of describing the map. This is also well known in
%the realm of ring isomorphisms between spaces of continuous
%functions: assuming that $X$ and $Y$ are realcompact and $T:C(X)
%\ra C(Y)$ is a ring isomorphism, there exists a homeomorphism
%$h:Y \ra X$ such that \[(Tf)(y)= (f \circ h)(y)\] for every $f \in
%C(X)$ and $y \in Y$. In
%the same way, a close description can be given for {\em linear} biseparating
%maps: if $T: C(X) \ra C(Y)$ is such a map, then there
%exist a homeomorphism $h:Y \ra X$ and a function $a \in C(Y)$
%with $a(y) \neq 0$ for every $y \in Y$ such that \[(Tf) (y) =
%a(y) (f \circ h) (y)\] for every $f \in C(X)$ and $y \in Y$ (\cite{ABN1, J}). 

        Of course, if we study ring isomorphisms between spaces
$C^* (X)$ and $C^* (Y)$, the conclusions we obtain are in
general poorer, in the sense that we can no longer conclude the
existence of a homeomorphism between the realcompactifications of $X$ and $Y$, but
just that the Stone-\v{C}ech compactifications of $X$ and $Y$
are homeomorphic. And this because every function $f \in C^*
(X)$ can be extended to a continuous function in $C(\beta X)$.
This implies that every ring isomorphism between $C^* (X)$ and
$C^* (Y)$ can be regarded as one between $C(\beta X)$ and
$C(\beta Y)$. Similar arguments also apply to biseparating maps
between spaces of bounded continuous functions: the
characteristics of these spaces yield a poor link between $X$
and $Y$. In particular, as an example, if we take a realcompact space $X$ which is not compact and define $Y:=\beta X$, then we have that the map $T: C^* (X) \ra C^* (Y)$ sending each $f \in C^* (X)$ into its extension to $\beta X$ is clearly a biseparating map. Nevertheless $X$ and $\beta X$ are not homeomorphic. This example illustrates also the fact that, in general, the existence of a biseparating map from $C^* (X)$ onto $C^* (Y)$ does not imply the existence of a similar map between $C(X)$ and $C(Y)$. Indeed, if such a map existed, we would conclude in our example that $X$ and $\beta X$, which are both realcompact, should be homeomorphic, against our assumptions.

        And what about spaces of vector-valued functions? Of
course, in this context ring isomorphism does not make sense,
but a close concept as biseparating map can be introduced here.

        In particular, a natural question arises: what happens if we consider
biseparating maps $\bxy$? Clearly if $E$ and $F$ are
infinite-dimensional, we can no longer act as in the
scalar-valued case, because since $E$ and $F$ are not locally
compact, bounded continuous functions cannot be extended to
functions on
$\beta X$ or $\beta Y$ attaining values in $E$ or $F$. On the other hand, notice  that classical techniques involving the study of ideals are no longer useful because  our spaces  do not have a ring structure.

        What we obtain in this paper is a different conclusion.
In our new context we are in a position to assure that the
existence of a biseparating map between $\xe$ and $\yf$ allows us
to claim that not only the Stone-\v{C}ech compactifications of
$X$ and $Y$ are homeomorphic, but also their
realcompactifications, a result that, as we mentioned above, is far from being true in
the scalar case. On the other hand, this result is the best one we can expect: that is, we cannot expect a homeomorphism between both spaces $X$ and $Y$.
%But we can say even more for realcompact spaces $X$ and $Y$: when the biseparating %map is
%linear it is also 
%possible to give a complete description of it, as it happened in
%the scalar context. 

Of course, we do not restrict our study to these special spaces, but we also obtain similar results for a large class of spaces, many of them containing unbounded functions. For instance, if we deal with spaces which do not contain "many" functions, the conclusions  we get apply not to the spaces $X$ and $Y$, but to some compactifications of them. But even in these cases, the structure of $X$ and $Y$ provides sometimes a strong link between them. In particular, this happens when we study biseparating maps between spaces of uniformly continuous bounded functions.

%$Finally we would like to mention some other
%recent works which also have some relation to this kind
%of problems, as \cite{AAK, A2, A1} and \cite{HBN}.

\section{Definitions and notation}

All over the paper $X$ and $Y$ will be completely regular topological spaces, and $E$ and $F$ will be ${\Bbb K}$-normed spaces. $C(Y, F)$ and $C^{*} (Y, F)$ are defined, with the natural modifications, in the same way as $C(X,E)$ and $C^{*} (X,E)$ were defined at the beginning of this section. The same comments apply to $C(Y)$ and $C^{*} (Y)$. Finally, if $X$ and $Y$ are also complete metric spaces, we introduce  $C_u^{*} (X, E)$ and $C_u^{*} (X, E)$ as the spaces  of uniformly continuous bounded functions defined on $X$ and $Y$, respectively, and taking values in $E$ and $F$, respectively. Also in this case $C_u^{*} (X) = C_u^{*} (X, {\Bbb K})$, and $C_u^{*} (Y) = C_u^{*} (Y, {\Bbb K})$.

\begin{definition}
{\em Given $f \in C(X,E)$, we define the cozero set of $f$ as \[ c(f) := \{ x \in X : f(x) \neq 0\}.\]}
\end{definition}

For the following two definitions, we assume that
$A, B$ are subrings of $C(X)$ and $C(Y)$, respectively, and that $\bx \subset C(X,E)$, $\ay \subset C(Y,F)$ are an $A$-module and a $B$-module, respectively.

\begin{definition}
{\em A map $\txy$ is said to be {\em separating} if it is additive and
$c(Tf) \cap c(Tg) = \emptyset$ whenever $f, g \in \bx$ satisfy
$c(f) \cap c(g) = \emptyset$. Besides $T$ is said to be
{\em biseparating} if it is bijective and both $T$ and $T^{-1}$ are separating.} 
\end{definition}

Equivalently, we see that an additive map $\txy$ is separating if $\vc (Tf) (y) \vd \vc (Tg) (y) \vd =0$ for all $y \in Y$ whenever $f, g \in \bx$ satisfy $\vc f(x) \vd \vc g(x) \vd=0$ for all $x \in X$. Notice then that in particular every ring isomorphism between $C(X)$ and $C(Y)$ is clearly a biseparating map. Linearity of maps will be assumed at no point of this paper.

\begin{definition}
{\em Let $\txy$ be a map and suppose that $\gamma X$ is a compactification of $X$. A point $x \in
\gamma X$ is said to be a
support point of $y \in Y$ if, for every neighborhood $U$
of $x$ in $\gamma X$, there exists $f \in \bx$ satisfying $c(f) \subset U$
such that $ (Tf ) (y) \neq 0$.}
\end{definition}
       
For a continuous map $f: X \ra {\Bbb K}$,
% $f^{\upsilon X}: \upsilon X \ra {\Bbb K}$ and 
$f^{\beta X}: \beta X \ra {\Bbb K} \cup \{\infty\}$ stands for the continuous extension to %$\upsilon X$ (the realcompactification of $X$) and 
$\beta X$ (the Stone-\v{C}ech compactification of $X$) into 
%${\Bbb K}$ and 
${\Bbb K} \cup  \{\infty \}$. In particular, given a continuous map $f : X \ra E$, $\vc f \vd^{\beta X}$ will be the continuous extension to $\beta X$ of $\vc \cdot \vd \circ f : X \ra {\Bbb K} \cup \{\infty\}$. In the same way, if $\gamma X$ is a compactification of $X$, and $f: X \ra {\Bbb K}$ is a continuous function which can be continuously extended to a map from $\gamma X$ into ${\Bbb K} \cup \{\infty\}$, we will denote by $f^{\gamma X}$ this natural extension.
Also, for ${\bf e} \in E$, $\widehat{{\bf e}}$ will stand for the function constantly equal to ${\bf e}$. Finally, given $C \subset X$ and $D \subset \gamma X$, $\cl C$ and $\cgx D$ will be their closures in $X$ and $\gamma X$, respectively.

Assume that $A$ is a subring of $C(X)$ which separates each point of $X$ from each point of $\beta X$. In $\beta X$, we introduce the equivalence relation $\sim$, defined as $ x \sim y$ whenever $f^{\beta X} (x) = f^{\beta X} (y)$ for every $f \in A$. In this way, we obtain the quotient space $\gamma X:= \beta X / \sim$. It is easy to see that $\gamma X$ is a compactification of $X$, and that every $f \in A$ can be continuously extended to a map from $\gamma X$ into ${\Bbb K} \cup \{\infty\}$. On the other hand, this extension will be bounded if $f$ is bounded.

Suppose that $\bx \subset C(X,E)$ is an $A$-module, where $A$ is a subring of $X$ which separates each point of $X$ from each point of $\beta X$. We say that $\bx$ is {\em  compatible} with $A$ if, for every $x \in X$, there exists $f \in \bx$  with $f(x) \neq 0$, and if, given any points $x, y \in \beta X$ with $x \sim y$, we have $\vc f \vd^{\beta X} (x) = \vc f \vd^{\beta X} (y)$ for every $f \in \bx$. It is clear that, in this case,  for each $f \in \bx$, there exists a continuous extension $\vc f \vd^{\gamma X} : \gamma X \ra {\Bbb K} \cup \{\infty\}$  of $\vc \cdot \vd \circ f$ to the whole space $\gamma X$.

A subring $A \subset C(X)$ is said to be {\em strongly regular} if given $x_0 \in \gamma X$ and a nonempty closed subset $K $ of $\gamma X$ which does not contain $x_0$, there exists $f \in A$ such that $f^{\gamma X} \equiv 1$ on a neighborhood of $x_0$ and $f^{\gamma X} (K) \equiv 0$.

\bigskip

\noindent
{\bf Examples.}
1. Suppose that $X$ is a metric space. Then the spaces $C(X, E)$ and  $C^{*} (X, E)$ are both $C^{*} (X)$-modules compatible with $C^{*} (X)$. Also, in this case, it is easy to see that $\gamma X = \beta X$.

2. Now suppose that $X$ is a complete metric space. It is easy to see that $C_u^{*} (X, E)$  is a $C_u^{*} (X)$-module compatible with $C_u^{*} (X)$. In this case,
in general, $\gamma X \neq \beta X$. On the other hand, it is immediate that we can embed isometrically our space $ C_u^{*} (X) $ in $C( \gamma X)$.  Also $ C_u^{*} (X) $ is a closed subalgebra of $C(\gamma X)$ which separates points, it contains constants and, when ${\Bbb K} ={\Bbb C}$, it is closed under complex conjugation. Then, by the Stone-Weierstrass theorem, it coincides with $C(\gamma X)$. Consequently $ C_u^{*} (X) $  is a strongly regular ring.

3. Suppose that $\Omega$ is a (not necessarily bounded) open subset of ${\Bbb R}^p$ ($p \in {\Bbb N}$). For a Banach space $E$, consider the set  $\coe$ of all $E$-valued functions whose partial derivatives up to the order $n$ exist and are continuous  ($n \ge 1$). Take $A= C^n (\Omega, {\Bbb R})$, which is a strongly regular ring  (see for instance \cite[Corollary 1.2]{W}). On the other hand, it is straightforward to see that $\coe$ is an $A$-module and that, in this case $\gamma \Omega = \beta \Omega$.

4. In a similar way, if $\Omega$ is a {\em bounded} open subset of
${\Bbb R}^p$ ($p \in {\Bbb N}$) and $E$ is a Banach space, we define $C^n
(\bar{\Omega}, E)$ as the
subspace of $\coe$ of those functions
whose partial derivatives up to order $n$
admit continuous extension to the boundary of $\Omega$. Clearly, we can view $C^n
(\bar{\Omega}, E)$  as a subset of $C(\clp \Omega, E)$
As above,  $A= C^n (\bar{\Omega}, {\Bbb R})$ is a strongly regular ring for which $\gamma \Omega = \clp \Omega$, and $C^n
(\bar{\Omega}, E)$ is an $A$-module.

\bigskip

A subring $A \subset C(X)$ is said to be {\em normal} if given two disjoint closed subsets $K, L$ of $\beta X$, there exists $f \in A$, $0 \le f \le 1$, satisfying $f^{\beta X}(K) \equiv 1$ and $f^{\beta X} (L) \equiv 0$. It is said to be {\em local} if $f \in C(X)$ belongs to $A$ whenever for every $x \in X$ there exist an open neighborhood $U(x)$ of $x$ and $f_x \in A$ such that $f \equiv f_x$ in $U(x)$.

	It is clear that each normal subring $A$ of $C(X)$ is also strongly regular, and that every $A$-module $\bx \subset C(X,E)$ must be compatible with $A$.

        In this paper we will assume we are in one of the following three situations.

\begin{itemize}
\item {\bf Situation 1.} $\bx \subset C(X,E)$ and $\ay \subset C(Y,F)$ are an
$A$-module and a $B$-module compatible with $A$ and $B$, respectively, where $A \subset
C(X)$ and $B \subset C(Y)$ are strongly regular rings. Also,  in the case when $\gamma X \neq \beta X$ and $\gamma Y \neq \beta Y$, we also assume that for every $x \in \beta X$ and $ y \in \beta Y$, there exist $f \in \bx$ and $g \in \ay$ satifying $\vc f \vd^{\beta X} (x) \neq 0$, $\vc g \vd^{\beta Y} (y) \neq 0$.
\item {\bf Situation 2.} $\bx \subset C(X,E)$ is an $A$-module and $\ay \subset C(Y,F)$ is a $B$-module, where $A \subset C(X)$ and $B \subset C(Y)$ are normal
local rings. 
\item {\bf Situation 3.} $E$ and $F$ are infinite-dimensional, and $\bx= C^* (X,E)$, $\ay = C^* (Y,F)$.
\end{itemize}

Notice that when we are in Situations 2 or 3, then we are also in Situation 1. Also it is clear that when we are in Situations 2 and 3, then $\gamma X =\beta X$ and $\gamma Y = \beta Y$.

We will denote by $\upsilon X$ and $\upsilon Y$ the realcompactifications of $X$ and $Y$, respectively.

	All over the paper the word "homeomorphism" will be synonymous with "surjective homeomorphism".

\section{Main results}

We first state a general result. Even if we may not have many functions in our spaces we can link the structures of some compactifications of $X$ and $Y$.

\begin{theorem}\label{3}
Suppose that we are in Situation 1. If $\txy$ is a  biseparating map, then $ \gamma X$ and
$\gamma Y$ are homeomorphic.
\end{theorem}

In some contexts, such as when $X$ and $Y$ have special structures, we can even ensure the existence of a homeomorphism between both spaces.

\begin{corollary}\label{txa}
Assume that we are in Situation 1. Suppose that, for $\bx$ and $\ay$, we have $\gamma X = \beta X$, and $\gamma Y = \beta Y$, respectively. If $X$ and $Y$ are first-countable spaces and $\txy$ is  biseparating,
then $X$ and $Y$ are homeomorphic. In particular, if $X$ and $Y$
are open subsets of ${\Bbb R}^p$ and ${\Bbb R}^q$, respectively,
then $p=q$.
\end{corollary}

\begin{corollary}\label{unif}
Let $\bx = C_u^{*} (X, E)$, $\ay = C_u^{*} (Y, F)$. If $\txy$ is a biseparating map, then $X$ and $Y$ are uniformly homeomorphic, that is, there exists a homeomorphism $h: Y \ra X$ such that both $h$ and $h^{-1}$ are uniform maps.
\end{corollary}

Notice that Corollary~\ref{txa} applies to Examples 1, 3 and 4. In a similar way, Corollary~\ref{unif} applies to Example 2.

Finally we state the main result of the paper, which applies to a large family of spaces of vector-valued functions.

\begin{theorem}\label{5}
Suppose that we are in Situations 2 or 3. If $\txy$ is  biseparating, then $\upsilon X$
and $\upsilon Y$ are homeomorphic.
\end{theorem}

In general we cannot conclude, in Theorem~\ref{5}, a statement like "$X$ and $Y$ are homeomorphic", 
as the following example shows.

\bigskip

\noindent
{\bf Example.}
Take any space $X$ which is not realcompact (for instance $W( \omega_1):= \{\sigma : \sigma < \omega_1\}$, where $\omega_1$ denotes the first uncountable ordinal; see \cite[5.12]{GJ}), and a realcompact normed space $E$. Recall that every normed space (or, more generally, every metrizable space) of nonmeasurable cardinal is realcompact (\cite[p. 232]{GJ}), and the assumption that all cardinal numbers are nonmeasurable is consistent with the axioms of set theory (\cite[p. 217]{En}). Next, each continuous map $f: X \ra E$ can be extended to a continuous map $f^{\upsilon X} : \upsilon X \ra E$. Clearly the map sending each $f \in C(X,E)$ into $f^{\upsilon X} \in C( \upsilon X, E)$ is biseparating but $X$ and $\upsilon X$ are not homeomorphic.

\section{Proofs I}

        In this section we assume that we are in Situation 1.

\begin{lemma}\label{a}
For any $x \in  X$, if $U$ is an open neighborhood of $x$ in
$\gamma X$, then
there exists $f \in \bx$ such that $x \in c(f)$ and $c(f) \subset U$.
\end{lemma}

{\em Proof.} 
Take $g \in \bx$ such that $x \in c(g)$, and $k
\in A$ such that $k(x) =1$ and
$k^{\gamma X} \equiv 0$ outside $U$. It is easy to see that $f:=
gk \in \bx$ does the job.
\hfill $\Box$

\begin{lemma}\label{A1}
Let $\txy$ be a  biseparating map. For $f, g \in \bx$, if $ c(f)
\subset c(g)$, then $c(Tf) \subset \cgy c(Tg)$.
\end{lemma}

{\em Proof.} 
Suppose that we can find $f, g \in \bx$ with $ c(f)
\subset c(g)$ and $c(Tf) \not\subset \cgy c(Tg)$. Take $y \in
c(Tf)$, $y \notin \cgy c(Tg)$. By Lemma~\ref{a} there exists
$k \in \ay$  such that $k(y) \neq 0$ and  $ c(k) \cap  \cgy
c(Tg) = \emptyset$. Consequently  $ c(k)
\cap  c(Tg) = \emptyset$ and, since $T^{-1}$ is  separating, we deduce
that $c(T^{-1} k) \cap c(g) = \emptyset$. Using the fact that
$T$ is  separating we deduce that $c(k) \cap c(Tf) =
\emptyset$.
Since this is a contradiction, we conclude that $c(f)
\subset \cgy c(Tg)$.
\hfill $\Box$

\begin{lemma}\label{1}
Let $\txy$ be a 
biseparating map. Then for each $y \in  Y$, there exists a unique
support point of $y$ in $ \gamma X$.
\end{lemma}

{\em Proof.} 
Take any $y \in Y$ and define $I_y := \{ f \in \ay : y \in
c( f  )\}$. Now consider $H(y) := \bigcap_{f \in
I_y} \cgx c(T^{-1} f)$.

\medskip

        {\bf Claim 1.} {\em $H(y)$ is nonempty.}

        We are going to see that the family $\{\cgx c(T^{-1} f)
: f \in I_y\} $ satisfies the 
finite intersection property. Take $f_1, f_2 , 
\ldots f_n \in I_y$. By Lemma~\ref{a}, we can consider $f 
\in I_y$ such that $ c( f) \subset \bigcap_{i=1}^n c(\vc f_i \vd^{\gamma Y})$.
This gives us $ c( f) \subset \bigcap_{i=1}^n c(f_i )$.
Then, by Lemma~\ref{A1} applied to $T^{-1}$, we have that
$c(T^{-1} f) \subset \cgx c(T^{-1} f_i)$ for every $i \in \{1,
2, \ldots, n\}$,
which implies that $\bigcap_{i=1}^n \cgx c(Tf_i) \neq
\emptyset$.  By the compactness of $\gamma X$,
we have that $H(y)$ is nonempty, and the claim is proved.

\medskip

        {\bf Claim 2.} {\em For each $y_0 \in Y$, $H(y_0)$ consists of
just one point.}

        Assume on the contrary that there exist two different
points $x_1 ,x_2 \in H(y_0)$. Take $f_0 \in \ay$ such that $f_0
(y_0) \neq 0$. By the definition of $H(y_0)$, we have that both
$x_1$ 
and $x_2$ belong to $\cgx c(T^{-1} f_0)$. Now consider a closed neighborhood $U_2$ of
$x_2$ such that $x_1 \notin U_2 $. Next, since $A$ is strongly regular, we can take $g_1
\in A$ 
such that $g_1^{\gamma X} \equiv 1$ on a neighborhood $U_1$ of
$x_1$, and
$g_1^{\gamma X} \equiv 0$ on $ U_2$. It is
clear that \[f_0 = T(g_1 (T^{-1} f_0)) + T( (1-g_1 ) (T^{-1} f_0)),\]and
consequently $y_0$ belongs to \[ c(T(g_1 (T^{-1} f_0)))\] or to
\[c(T((1-g_1 )(T^{-1} f_0))).\] We assume without loss of generality
that $y_0 \in  c(T(g_1 (T^{-1} f_0)))$.

        Now, since $y_0 \in  c(T(g_1 (T^{-1}
f_0)))$ and $x_2 \in H(y_0)$, then $x_2 $ belongs to $\cgx c(g_1
(T^{-1} f_0))$, which is not true by construction.

        We conclude that $H(y_0)$ contains just one point.

\medskip

        {\bf Claim 3.} {\em Given $y_0 \in Y$, if $x \in H(y_0)$, then
$x$ is a support point 
of $y_0$.}

        Consider an open neighborhood $U$ of $x \in H(y_0)$ in
$\gamma X$. We have to prove that there exists $g_0 \in \bx$ such
that $c(g_0) \subset U$ and $(Tg_0) (y_0) \neq 0$.

        Take $f_0 \in \ay$ such that $f_0 (y_0) \neq 0$.  Of
course,  if $c(T^{-1} f_0)$ is
contained in $U$, we get the result by defining $g_0 := 
T^{-1} f_0$, so we suppose this is not the case. 

        Take $g_1 \in A$ such that $\cgx c(g_1)
\subset U$ and $g_1^{\gamma X} \equiv 1$ on a neighborhood of
$x$. It is clear that \[T^{-1} f_0 = g_1 T^{-1} f_0 + (1 - g_1)
T^{-1} f_0,\] and, as above, as a consequence we have that \[(T (g_1 T^{-1}
f_0)) (y_0) \neq 0\] or  \[(T ((1- g_1) T^{-1}
f_0)) (y_0) \neq 0.\] But notice that if the latter holds, since
$x \in H(y_0)$ we should have $x \in \cgx c((1-g_1) T^{-1}
f_0)$, which is not the case. Consequently, defining $g_0 := g_1
T^{-1} f_0$ we are done.
\hfill $\Box$ 

        The previous lemma allows us to define a map $h: Y \ra
\gamma X$, sending
each point $y \in Y$ into its support point $ h (y) \in \gamma
X$.

\begin{lemma}\label{13}
Let $\txy$ be a  biseparating map. Suppose that $h(y)
= x$ for some $y \in Y$, and that $f \in \bx$
satisfies $\vc f \vd^{\gamma X} \equiv 0$ on a neighborhood of
$x$. Then $Tf \equiv 0$ on a neighborhood of $y$.
\end{lemma}

{\em Proof.} 
Take an open neighborhood $U$ of $x$ in $\gamma X$
such that $\vc f \vd^{\gamma X} \equiv 0 $ in $U$. By the
definition of support point, we can take $g \in \bx$ such
that $ c(g) \subset U$ and $(Tg) (y) \neq
0$. Since $T$ is  biseparating, and $ c(f) \cap  c(g) = \emptyset$,
we deduce that $Tf \equiv 0$ in $c(Tg)$, which is a neighborhood of $y$.
\hfill $\Box$

        Now the following corollary follows easily.

\begin{corollary}\label{sat}
Let $\txy$ be a  biseparating map. Suppose that $y \in Y$, and that $f
\in \bx$ satisfies $ (Tf) (y) \neq 0$. Then $h (y) \in
\cgx c(f)$.
\end{corollary}

\begin{lemma}\label{A2}
Given a  biseparating map $\txy$, the associated map $h:Y \ra
\gamma X$ is continuous, and its range is dense in $\gamma X$.
\end{lemma}

{\em Proof.} 
We shall see first that $h$ is continuous at every point of $Y$.
Take $y_0 \in Y$ and an open neighborhood $U$ of $h(y_0)$ in
$\gamma X$. By the definition of support point, there exists $g
\in \bx$ such that $\cgx c(g) \subset 
U$ and $(Tg) (y_0) \neq 0$. Thus $c(Tg)$ is an open neighborhood
of $y_0$. By Corollary~\ref{sat}, $h (c(Tg))
\subset U$ and we are done.

        On the other hand, taking into account that $T$ is
injective, a similar reasoning yields that for every open subset
$U$ of $\gamma X$, there are points $h(y)$ in $U$, and
consequently the range of $h$ is dense in $\gamma X$.
\hfill $\Box$

\begin{lemma}\label{ext}
The map $h$ can be extended to a continuous map from $\gamma Y$ onto $\gamma X$.
\end{lemma}

{\em Proof.} 
Obviously, as a consequence of the previous lemma, we can extend $h$
to a continuous map $\widehat{h}$ from $\beta Y$ onto $\gamma X$, so the result is true if $\gamma Y = \beta Y$.

Thus we assume that $\gamma Y \neq \beta Y$. 

{\bf Claim 1.} {\em Let $y_0 \in \beta Y$. Given an open neighborhood $U$ of $\widehat{h} (y_0) $ in $\gamma X$, there exists $f \in \bx$ such that $c(f) \subset U$ and $\vc Tf \vd^{\beta Y} (y_0) \neq 0$.}

Notice that, by hypothesis, since $\gamma Y \neq \beta Y$, we have that there exists $g \in \ay$ such that $\vc g \vd^{\beta Y} (y_0) \neq 0$. Now take $k \in A$ such that $k^{\gamma Y} \equiv 0$ outside $U$ and $k^{\gamma Y} \equiv 1$ on an open  neighborhood, say $V$, of $\widehat{h}(y_0)$. It is clear that, if we define $f:= k T^{-1} g$, then $c(f) \subset U$. Also, $\widehat{h}^{-1} (V)$ is a neighborhood of $y_0$. Now, if we take $y \in \widehat{h}^{-1} (V) \cap Y$, then $h(y) \in V$, and consequently $\vc f- T^{-1} g \vd^{\gamma X} \equiv 0$ on a neighborhood of $h(y)$. By Lemma~\ref{13}, $(T f) (y) = g(y)$. Since this holds for every $y \in \widehat{h}^{-1} (V) \cap Y$, we conclude that $\vc Tf \vd^{\beta Y} (y_0 ) = \vc g \vd^{\beta Y} (y_0) \neq 0$.

\medskip

{\bf Claim 2.} {\em Given $y_1 , y_2 \in \beta Y$, if $y_1 \sim y_2$, then $\widehat{h} (y_1) = \widehat{h} (y_2)$.}

Suppose that $y_1 \sim y_2$, $y_1 \neq y_2$, and that $\widehat{h} (y_1) \neq \widehat{h}(y_2)$. Take disjoint open subsets $U$ and $V$ of $\widehat{h} (y_1) $ and $\widehat{h} (y_2) $, respectively, in $\gamma X$. Now, by Claim 1, there exist $f_1 , f_2 \in \bx$ such that $c(f_1) \subset U$, $c(f_2) \subset V$, and $\vc Tf_1 \vd^{\beta Y} (y_1) \neq 0$,  $\vc Tf_2 \vd^{\beta Y} (y_2) \neq 0$. But, since $T$ is separating and $c(f_1) \cap c(f_2) = \emptyset$ , we deduce that $\vc (Tf_1) (y) \vd \vc (Tf_2) (y) \vd =0$ for every $y \in Y$. Clearly this must force to $\vc Tf_1 \vd^{\beta Y} (y_2) =0$. This contradicts the fact that $\ay$ is compatible with $B$, which means in particular that $\vc g \vd^{\beta Y} (y_1) = \vc g \vd^{\beta Y} (y_2)$ for every $g \in \ay$.

\medskip

Finally, because of Claim 2, given $y \in \gamma Y$ we can define the image of $y$ as the image by $\widehat{h}$ of any of the elements of its equivalence class. It is clear that this determines a surjective continuous map from $\gamma Y$ onto $\gamma X$ which is an extension of $h$.
\hfill $\Box$ 

\bigskip

The extension map given in Lemma~\ref{ext} will also be called $h$.

\bigskip

{\em Proof of Theorem~\ref{3}.} 
We will show that $h$ is a homeomorphism. We are going to find an inverse for the map $h$. It is clear
that since $T$ is  biseparating, we can construct a function $k:
\gamma X \ra \gamma Y$ associated to $T^{-1}$, which is an extension
of a map from $X$ into $\gamma Y$ sending each point of $X$ into
its support point for $T^{-1}$. We just have to prove that $k$
is the inverse map of $h$. 

\medskip

        {\bf Claim 1.} {\em  If $y_0 \in Y$, then $k(h(y_0)) = y_0$.}

        Assuming the contrary, suppose that
$y_0 \in Y$ satisfies $y_0 \neq k(h(y_0))$. Then consider $U$,
$V$ two open neighborhoods in $\gamma Y$ of $y_0$ and $k(h(y_0))$
respectively, such that $U \cap V = \emptyset$.

Next, applying Lemma~\ref{a}, take $f
\in \ay$ such that $c(f) \subset U$ and $ f (y_0) \neq 0$. By
Corollary~\ref{sat}, we have that 
$h(y_0) $ belongs to $\cgx c(T^{-1} f)$.

        On the other hand, as $k: \gamma X \ra \gamma Y$ is continuous there
exists an open neighborhood $U_1$ of $h(y_0)$ in $\gamma X$ such that if
$x \in U_1$, then $k(x) \in V$. Also since $h(y_0) \in \cgx
c(T^{-1} f)$, then $U_1 \cap c(T^{-1} f) \neq \emptyset$, so we may
take \[x_0 \in U_1 \cap c(T^{-1} f).\]

        Since $x_0 \in U_1 $, $k(x_0) \in V$. But $x_0$ also belongs to $X$, and then, by the
definition of support point, there exists $g \in \ay$ 
such that $c(g) \subset V$ and $(T^{-1} g) (x_0) \neq 0$.

        So we have first that $c(f) \cap  c(g) = \emptyset$. But it
is clear that \[x_0 \in c(T^{-1} f) \cap c(T^{-1} g),\] contradicting the
fact that $T^{-1}$ is  separating.

        As a consequence Claim 1 is proved.

\medskip

        {\bf Claim 2.} {\em If $y_0 \in \gamma Y- Y$, then $k(h(y_0)) = y_0$.}

        Take a net $(y_{\alpha})$ in $Y$ converging to $y_0
\in \gamma Y- Y$. Since both $h$ and $k$ are continuous, we
have that the net $(k(h(y_{\alpha})))$ converges to $k(h(y_0))$.
But by Claim 1, $k(h(y_{\alpha})) = y_{\alpha}$ for every
$\alpha$. This implies that it converges to $y_0$, and
consequently $k(h(y_0)) = y_0$. So Claim 2 is proved.

\medskip

        As a consequence we easily conclude that $k$ is the inverse map
of $h$, and that both are homeomorphisms.
\hfill $\Box$ 

\bigskip

{\em Proof of Corollary~\ref{txa}.} 
By the previous theorem, we have that $\beta X$ and $\beta Y$
are homeomorphic. Also, since the only points of
$\beta X$ having a countable base of neighborhoods belong to $X$ (\cite[9.7]{GJ}), and the same applies to $Y$, we easily
conclude that $X$ and $Y$ are homeomorphic. Finally, in the
special case when
 $X$ and $Y$
are open subsets of ${\Bbb R}^p$ and ${\Bbb R}^q$, respectively,
  this fact implies that $p=q$ (see for instance \cite[p. 120]{Fe}). 
\hfill $\Box$

\bigskip

{\em Proof of Corollary~\ref{unif}.} 
In \cite[Lemma 3.4]{AF2}, it is proved, in a different context, that every point of $X$ is a $G_{\delta}$-set in $\gamma X$, and that, on the contrary, no point in $\gamma X - X$ is a $G_{\delta}$-set in $\gamma X$. Consequently it follows that $h$ is a homeomorphism from $Y$ onto $X$.

Then we have a map $S: C( \gamma X) \ra C (\gamma Y)$, defined as $(Sf) (y) = f (h(y))$ for every $f \in C(\gamma X)$ and every $y \in \gamma Y$. This map is easily seen to be bijective. Now, since $ C_u^{*} (X) $  and $C_u^{*} (Y)$ can  be identified with $C(\gamma X)$ and $C(\gamma Y)$, respectively (see Example 2), we may consider $S$ as a map from $C_u^{*} (X)$ onto $C_u^{*} (Y)$, where it is  also defined as $Sf = f \circ h$. As a consequence,
for each $f \in  C_u^{*} (X) $, $f \circ h$ belongs to  $C_u^{*} (Y) $.

But, on the other hand, we can prove as in \cite[Theorem 2.3]{LL} (see also the Remark after it) that if $f \circ h \in C_u^{*} (Y)$ whenever $f \in C_u^{*} (X)$, then $h$ is uniformly continuous. Since the same process works also for $h^{-1}$, then the theorem is proved.
\hfill $\Box$

\section{Proofs II}

        In this section we will assume that we are in Situations 2 or 3.
% (stated here just for the space $\bx$):

Recall that in this context,
if we look at the equivalence relation introduced in Section 2, then $\gamma X = \beta X$ and $\gamma Y = \beta Y$, so previous results can be applied here.

We start with a result whose proof is easy from
Lemma~\ref{13} and the fact that  $h$ is a homeomorphism.

\begin{lemma}\label{cc}
Let $\txy$ be a  biseparating map. Suppose that $y \in \beta Y$,
and that $f \in \bx$ 
satisfies $\vc f \vd^{\beta X} \equiv 0$ on an open set $U$ of
$\beta X$. Then $\vc Tf \vd^{\beta Y} \equiv 0$ on the open set $h^{-1} (U) \subset \beta Y$.
\end{lemma}

\begin{prop}\label{2}
If $\txy $ is a  biseparating  map, then for every $y \in Y$,
$h(y) \in \upsilon X$.
\end{prop}

{\em Proof.} 
Take $y_0
\in  Y$, and
suppose that $h (y_0) \in \beta X 
- \upsilon X$. Then there exists
a sequence $(U_n)$ of open neighborhoods of $h(y_0)$ in $\beta X$
such that $\clx U_{n+1} \subset U_n$, $\clx U_{n+1} \neq U_n$,
for every $n \in {\Bbb N}$ 
and $X \cap \bigcap_{n=1}^{\infty} U_n = \emptyset$ (see for
instance \cite[Theorem 3.11.10]{En}). It is clear that \[h(y_0)
\in \clx \bigcup_{n \in {\Bbb N}} (U_{n} - \clx U_{n+2}).\]

        Since $h$ is a homeomorphism, we deduce that $y_0 \in
\cly \bigcup_{n \in {\Bbb N}} h^{-1} (U_{n} - \clx 
U_{n+2})$. Consequently,
if we define
\begin{eqnarray*}
V_1 &:=& 
\bigcup_{n \in {\Bbb N}} h^{-1} (U_{4n} - \clx U_{4n+2}),\\ V_2 &:=&
\bigcup_{n \in {\Bbb N}} h^{-1} (U_{4n-1} - \clx U_{4n+1}), \\ V_3 &:=&
\bigcup_{n \in {\Bbb N}} h^{-1} (U_{4n-2} - \clx U_{4n}),
\end{eqnarray*}
and \[V_4 :=
\bigcup_{n \in {\Bbb N}} h^{-1} (U_{4n-3} - \clx U_{4n-1}),\]
then we have that
$y_0 $ belongs to one of the sets $\cly V_1, \cly V_2 , \cly
V_3 $ or $\cly V_4$. We assume without loss of generality that
$y_0 \in \cly V_1$. 

At this point we split the proof into two cases.
\begin{itemize}
\item {\em Case 1. Assume we are in Situation 3.}

It is clear that, if for each $n \in
{\Bbb N}$, we define $W_n := h^{-1} (U_{4n} - \clx U_{4n+2})$, then $\cl W_n
\cap \cl W_m = \emptyset$, for $n \neq m$. 
Consider now a
sequence $(g_n)$ in $B$ such that, for every $n \in
{\Bbb N}$, $0 \le g_n \le 1$,
$g_n \equiv 1$ on $Y \cap W_n$, and $g_n \equiv 0$ outside $Y
\cap h^{-1} (U_{4n-1} - \clx U_{4n+3})$. 

Applying Riesz's
Lemma (\cite[Theorem 1.3.2]{Lr}), we can take a sequence 
$({\bf e}_n)$ of norm one points in $F$ 
such that $\vc {\bf e}_n - {\bf e}_m \vd \ge 1/2$ for $n \neq m$. 

        It is also clear that, since $T^{-1}$ is  separating, then
\[c(T^{-1} (g_n {\bf e}_n)) \cap c(T^{-1} (g_m
{\bf e}_m) )= \emptyset\] if $n 
\neq m$. Define $f :=
\sum_{n=1}^{\infty} T^{-1} (g_n {\bf e}_n)$. Notice that, for each $x \in X$, $f(x)$ belongs to $E$.

\medskip

        {\bf Claim 1.} {\em  $f$ is continuous.}

        Notice that each function $\vc g_n
{\bf e}_n \vd^{\beta Y} \equiv 0$ outside $\cly h^{-1} (U_{4n-1} - \clx
U_{4n+3})$, which
implies by Lemma~\ref{cc} applied to $T^{-1}$ that \[\vc T^{-1} (g_n
{\bf e}_n ) \vd^{\beta X} \equiv 0\] outside $\clx (U_{4n-1} - \clx
U_{4n+3})$. This is, we have that
for every $n \in 
{\Bbb N}$, $c(T^{-1} (g_n {\bf e}_n)) \subset \clx (U_{4n-1} -
\clx U_{4n+3})$. Now given any $x \in X$, there exists an open
neighborhood  $U$ of $x$ in $X$ such that there are just
a few numbers $k \in {\Bbb N}$ satisfying $U \cap \clx (U_{4k-1} - \clx
U_{4k+3}) \neq \emptyset$, due to the construction of the sequence
$(U_n)$. This proves Claim 1.

\medskip

        {\bf Claim 2.} {\em  $f$ belongs to $\bx$.}

        To prove it, we just need to show that $f$
is bounded. Suppose on the contrary that the sequence $( \vc T^{-1}
(g_n {\bf e}_n ) \vd )$ is not bounded. For each $n \in {\Bbb N}$,
set $a_n := \vc
T^{-1} (g_n {\bf e}_n ) \vd$. Since the sequence $(a_n)$ is not
bounded, we can extract a subsequence (which without loss of
generality we shall assume it to be the whole
$(a_n)$) with the property that $a_n \ge n^3$ for every $n \in
{\Bbb N}$. Next
consider the map \[g:= \sum_{n \in {\Bbb N}} \frac{g_n
{\bf e}_n}{n^2}.\]It is clear that, since each $\vc g_n
{\bf e}_n \vd \le 1$, then $g$ belongs to $\ay$. Consequently $T^{-1}
g$ exists and is a bounded
function on $X$. We are going to see that this is not true,
obtaining a contradiction. To this end take, for each $n \in
{\Bbb N}$, $x_n \in X$ such that \[\vc (T^{-1}
(g_n {\bf e}_n )) (x_n) \vd \ge \frac{a_n}{2} \ge
\frac{n^3}{2}.\] Now notice
that, for each $n_0 \in {\Bbb N}$, \[g = \frac{g_{n_0}
{\bf e}_{n_0}}{n_0^2} + \sum_{n \neq n_0} \frac{g_n
{\bf e}_n}{n^2},\] where \[ c \pa \frac{g_{n_0}
{\bf e}_{n_0}}{n_0^2} \rb \cap  c \pa \sum_{n \neq n_0} \frac{g_n
{\bf e}_n}{n^2} \rb = \emptyset.\] Since $T^{-1}$ is  separating,
we deduce that \[c \pa T^{-1} \pa \frac{g_{n_0}
{\bf e}_{n_0}}{n_0^2} \rb \rb \cap c\pa  T^{-1} \sum_{n \neq n_0} \frac{g_n
{\bf e}_n}{n^2} \rb = \emptyset,\]and this fact implies that
\[\vc (T^{-1} g) (x_{n_0}) \vd = \frac{1}{n_0^2} \vc (T^{-1} (g_{n_0}
{\bf e}_{n_0})) (x_{n_0}) \vd \ge \frac{n_0}{2}.\] Since this applies to
every $n_0 \in {\Bbb N}$, we deduce that $T^{-1} g$ is not
bounded. This contradiction allows us to conclude that $f$ is
bounded, and Claim 2 is proved.

\medskip

        Consequently we have that $Tf$ belongs to $\ay$,
and for each $n \in {\Bbb N}$, \[\vc f - T^{-1} (g_n
{\bf e}_n)\vd^{\beta X} \equiv 0\] on
$U_{4n} - \clx U_{4n+2}$, which implies by Lemma~\ref{cc} that
\[Tf \equiv g_n {\bf e}_n \equiv {\bf e}_n \] in $Y \cap W_n$.
Also, suppose that 
$(Tf) (y_0) = {\bf e}_0 \in F$. Then there exists a neighborhood $U$
of $y_0$ in $Y$ such that if $y, y' \in U$, then \[\vc (Tf) (y) - (Tf)
(y') \vd < \frac{1}{4}.\] 

        On the other hand, recall that $y_0
\in \cly V_1$ and, since $h(y_0) \notin \clx (U_n - \clx
U_{n+2})$ for any $n \in {\Bbb N}$, then $y_0 \notin \cly W_n$ for any
$n \in {\Bbb N}$. This implies that for every $k \in {\Bbb N}$,
\[y_0 \in \cly \bigcup_{n \ge k} W_n.\] As a consequence, there
are $n_1, n_2 \in 
{\Bbb N}$, $n_1 \neq n_2$, 
such that $U \cap W_{n_1} \neq \emptyset \neq U \cap W_{n_2}$, and
then there are $y_{n_1} \in U \cap W_{n_1}$, $y_{n_2} \in U \cap
W_{n_2}$. Also \[\frac{1}{4} > \vc (Tf) (y_{n_1})  - (Tf) (y_{n_2}) \vd = \vc
{\bf e}_{n_1} - {\bf e}_{n_2} \vd
\ge \frac{1}{2},\] which is impossible.

\item {\em Case 2. Assume we are in Situation 2}.

In this case we follow a similar pattern of proof to the given above.
Now consider a
sequence $(f_n)$ in $A$ such that, for every $n \in {\Bbb N}$,
$c(f_n) \subset X \cap U_{4n-2}$,
and $f_n (x) =1$ for every $x \in
X \cap U_{4n}$. Define $g := \sum_{n=1}^{\infty} f_n $.
Since $A$ is local, it is easy to check that $g$ belongs to $A$.
Also it is
easy to see that for every $n \in {\Bbb N}$, $g$ is
constantly equal to $n$ on
$X \cap (U_{4n} - \clx U_{4n+2})$. Next take $f \in \bx$ such
that $(Tf) (y_0) = {\bf f}_0 \neq 0$.
Suppose that $ T(gf)  (y_0) = {\bf f}_1 \in
F$. Consider $n_0 \in
{\Bbb N}$, $n_0 \vc {\bf f}_0 \vd /2 > \vc {\bf f}_1 \vd +1$, and an open
neighborhood $U(y_0)$ of $y_0$ in $\beta Y$ such that $h(U(y_0)) \subset
U_{4n_0} \cap V$, where \[V:= \left\{x \in \beta X : \vc T f
\vd^{\beta Y} (h^{-1} (x)) > \frac{{\bf f}_0}{2}.\right\}\] Since
$h(y_0 )$ belongs to $\clx V_1$, then as above for every $k \in {\Bbb N}$,
\[y_0 \in \cly \bigcup_{n \ge k} W_n.\] Now it is 
not difficult to see that there exists $k \in
{\Bbb N}$, $k \ge n_0$, such that $h(U(y_0)) \cap (U_{4k} - \clx
U_{4k+2})$ is nonempty. Then if for $y_1 \in U(y_0) \cap Y$, $h(y_1)$ belongs
to $U_{4k} - \clx 
U_{4k+2}$, we have that $\vc gf -kf \vd^{\beta X}$ is constantly
equal to zero in
a neighborhood of $h(y_1)$. This means, by Lemma~\ref{cc}, that $ T(gf -kf)
 (y_1) = 0$,
and consequently, \[\vc T (gf)
 (y_1) \vd = k \vc (Tf) (y_1) \vd \ge n_0 \frac{\vc {\bf f}_0
\vd}{2} \ge \vc
{\bf f}_1 \vd +1.\]

        Since this happens for every open neighborhood of $y_0$,
we deduce that $ T(gf) $ is not continuous, which is not
possible.

\end{itemize}

In both cases, we conclude that $h(y)$ belongs to $\upsilon X$ for
every $y \in Y$.
\hfill $\Box$ 

\bigskip

{\em Proof of Theorem~\ref{5}.} 
We have that the restriction of $h$ to
$\upsilon Y$ is continuous. Also, by Proposition~\ref{2}, $h(y)$
belongs to $\upsilon X$ for every $y \in Y$. Since $\upsilon X$
is realcompact, we deduce that $h(\upsilon Y)$ is contained in
$\upsilon X$. Since $h^{-1}$, for the same reason,
maps elements of $\upsilon X$ into elements of $\upsilon Y$, and
$h^{-1}$ is continuous, we conclude that the restriction of $h$
to $\upsilon Y$ is a homeomorphism onto $\upsilon X$.
\hfill $\Box$ 

\bigskip

The author wishes to thank the referee and Professor K. Jarosz for some suggestions which improved this paper.

\vspace{.15in}

        {\footnotesize {\sc Departamento de Matem\'aticas,
Estadistica y Computaci\'on, Facultad de Ciencias, Universidad
de Cantabria, Avenida de los Castros s.n., E-39071, Santander, Spain.}}

        {\footnotesize {\em E-mail address}: araujoj@@unican.es}

\end{document}